\documentclass[12pt]{article}

\newcommand{\blind}{1}

\usepackage{xr-hyper}
\usepackage{verbatim}
\usepackage{float}
\usepackage{rotfloat}
\usepackage{mathrsfs}
\usepackage{mathtools}
\usepackage{url}
\usepackage{enumitem}
\usepackage{amsmath}
\usepackage{amsthm}
\usepackage{amssymb}
\usepackage{graphicx}
\usepackage[authoryear]{natbib}
\usepackage[unicode=true,pdfusetitle,
 bookmarks=true,bookmarksnumbered=false,bookmarksopen=false,
 breaklinks=false,pdfborder={0 0 1},backref=false,colorlinks=false]
 {hyperref}
\hypersetup{
 colorlinks,allcolors=blue}

\makeatletter

\theoremstyle{plain}
\newtheorem{thm}{\protect\theoremname}[section]
\theoremstyle{plain}
\newtheorem{prop}{\protect\propositionname}[section]
\theoremstyle{plain}
\newtheorem{cor}{\protect\corollaryname}[section]
\theoremstyle{definition}
\newtheorem{defn}{\protect\definitionname}[section]
\theoremstyle{plain}

\@ifundefined{date}{}{\date{}}
\usepackage{xcolor}
\AtBeginDocument{%
\let\ref\autoref
\renewcommand\equationautorefname{\@gobble}

}

\makeatother

\providecommand{\corollaryname}{Corollary}
\providecommand{\definitionname}{Definition}
\providecommand{\lemmaname}{Lemma}
\providecommand{\propositionname}{Proposition}
\providecommand{\theoremname}{Theorem}

\global\long\def\ntoinf{n \rightarrow\infty}%
\global\long\def\toinf{\rightarrow\infty}%
\global\long\def\tozero{\rightarrow0}%
\global\long\def\sumjinf{\sum_{j=1}^{\infty}}%
\global\long\def\sumkinf{\sum_{k=1}^{\infty}}%
\global\long\def\sumkK{\sum_{k=1}^{K}}%
\global\long\def\sumin{\sum_{i=1}^{n}}%
\global\long\def\quadas{\quad a.s.}%
\global\long\def\bbR{\mathbb{R}}%
\global\long\def\bbN{\mathbb{N}}%
\global\long\def\bbH{\mathbb{H}}%
\global\long\def\bbE{\mathbb{E}}%
\global\long\def\bbF{\mathbb{F}}%
\global\long\def\hv{\hat{v}}%
\global\long\def\hmu{\hat{\mu}}%
\global\long\def\hV{\hat{V}}%
\global\long\def\hLambda{\hat{\Lambda}}%
\global\long\def\hlambda{\hat{\lambda}}%
\global\long\def\hphi{\hat{\phi}}%

\global\long\def\cA{\mathcal{A}}%
\global\long\def\cB{\mathcal{B}}%
\global\long\def\cC{\mathcal{C}}%
\global\long\def\cF{\mathcal{F}}%
\global\long\def\cG{\mathcal{G}}%
\global\long\def\cH{\mathcal{H}}%
\global\long\def\cI{\mathcal{I}}%
\global\long\def\cM{\mathcal{M}}%
\global\long\def\cN{\mathcal{N}}%
\global\long\def\cP{\mathcal{P}}%
\global\long\def\cS{\mathcal{S}}%
\global\long\def\cU{\mathcal{U}}%
\global\long\def\cX{\mathcal{X}}%
\global\long\def\sT{\mathscr{T}}%
\global\long\def\hsT{\mathscr{\hat{T}}}%
\global\long\def\cov{\operatorname{cov}}%
\global\long\def\pool{\text{pool}}%

\global\long\def\argmin{\operatornamewithlimits{arg\,min}}%

\global\long\def\conv#1{\mathrel{\mathop{\longrightarrow}\limits ^{#1}}}%

\global\long\def\sqrtn{\sqrt{n}}%
\global\long\def\oneovern{\frac{1}{n}}%
\global\long\def\Sinf{\cS^{\infty}}%
\global\long\def\ThmuS{T_{\hmu}\Sinf}%
\global\long\def\hcG{\hat{\cG}}%
\global\long\def\hcF{\hat{\cF}}%
\global\long\def\TmuS{T_{\mu}\Sinf}%
\global\long\def\logmu{\log_{\mu}}%
\global\long\def\loghmu{\log_{\hmu}}%
\global\long\def\pinSInf{p \in\Sinf}%
\global\long\def\norm#1{\left\lVert #1\right\rVert }%
\global\long\def\normE#1{\left\lVert #1\right\rVert _{\bbE}}%
\global\long\def\Norm#1{{\left\vert \kern-0.25ex  \left\vert \kern-0.25ex  \left\vert #1 \right\vert \kern-0.25ex  \right\vert \kern-0.25ex  \right\vert }}%

\global\long\def\rhotau{\rho_{\tau}}%
\global\long\def\mutau{\mu_{\tau}}%
\global\long\def\hmutau{\hmu_{\tau}}%
\global\long\def\psitau{\psi_{\tau}}%

\global\long\def\inner#1#2{\langle#1,#2\rangle}%
\global\long\def\innerE#1#2{\langle#1,#2\rangle_{\bbE}}%
\global\long\def\HS{\text{HS}}%
\global\long\def\innerp#1#2{\langle#1,#2\rangle_{p}}%
\global\long\def\innerHS#1#2{\langle#1,#2\rangle_{\HS}}%

\global\long\def\Dens{\text{Dens}}%
\DeclarePairedDelimiterX{\Inner}[2]{\langle\langle}{\rangle\rangle}{#1, #2}

\def\references{\bibliography{/Users/xdai/GoogleDrive/Papers/zotero}}

\begin{document}

\def\spacingset#1{\renewcommand{\baselinestretch}%
{#1}\small\normalsize} \spacingset{1}

\if1\blind
{
  \title{\bf Statistical Inference on the Hilbert Sphere with Application to Random Densities}
  \author{Xiongtao Dai\thanks{
    \texttt{xdai@iastate.edu}}\hspace{.2cm}\\
    Department of Statistics, Iowa State University, Ames, Iowa 50011 USA \\
    }
  \maketitle
} \fi

\if0\blind
{
  \bigskip
  \bigskip
  \bigskip
  \begin{center}
    {\LARGE\bf Statistical Inference on the Hilbert Sphere with Application to Random Densities}
\end{center}
  \medskip
} \fi

\bigskip
\begin{abstract}
The infinite-dimensional Hilbert sphere $\Sinf$ has been widely employed to model density functions and shapes, extending the finite-dimensional counterpart. 
We consider the Fréchet mean as an intrinsic summary of the central tendency of data lying on $\Sinf$. 
To break a path for sound statistical inference, we derive properties of the Fréchet mean on $\Sinf$ by establishing its existence and uniqueness as well as a root-$n$ central limit theorem (CLT) for the sample version, overcoming obstructions from infinite-dimensionality and lack of compactness on $\Sinf$. 
Intrinsic CLTs for the estimated tangent vectors and covariance operator are also obtained. 
Asymptotic and bootstrap hypothesis tests for the Fréchet mean based on projection and norm are then proposed and are shown to be consistent. 
The proposed two-sample tests are applied to make inference for daily taxi demand patterns over Manhattan modeled as densities, of which the square roots are analyzed on the Hilbert sphere.
Numerical properties of the proposed hypothesis tests which utilize the spherical geometry are studied in the real data application and simulations, where we demonstrate that the tests based on the intrinsic geometry compare favorably to those based on an extrinsic or flat geometry.
\end{abstract}

\noindent%
\vfill

\newpage
\spacingset{1.45} %

\section{Introduction\label{s:intro}}

We aim to develop statistical theory and inferential methods for analyzing a sample of random elements taking values on the Hilbert sphere $\Sinf=\{f\in\bbH\mid\norm f=1\}$, where $\bbH$ is a separable infinite-dimensional Hilbert space with inner product $\inner{\cdot}{\cdot}$ and norm $\norm{\cdot}$.
The spherical Hilbert geometry gives rise to invariance properties and efficient calculation of geometric quantities, for instance the geodesics, which promotes much of its prevailing applications to model different data objects:

\begin{itemize}[labelindent=0pt,wide] %
\item 
\emph{Probability distributions} supported on a compact interval $\cN$ is represented by density functions in $\Dens(\cN)=\{y : \cN \rightarrow \bbR \mid\int_{\cN}y(s)\,ds=1,\ y(s)>0,\, s\in\cN\}$. 
To analyze the densities in a reparametrization-invariant fashion, \citet{sriv:07} proposed a nonparametric Fisher--Rao metric with a reparametrization-invariant property, generalizing the well-known finite-dimensional version \citep{rao:45}.
The Fisher--Rao metric is conveniently expressed and calculated through representing the densities by the corresponding square root densities in $\cX=\{x : \cN \rightarrow \bbR \mid x(s)=\sqrt{y(s)},\,s\in\cN,\,y\in\Dens(\cN)\}$, where $\cX$ is the positive orthant of the Hilbert sphere $\Sinf$ modeled in $\bbH=L^{2}(\cN)$.
The square root density framework has been widely applied in modeling the orientation distribution functions in high angular resolution diffusion images \citep{chen:09-1,du:14} and time-warping functions \citep{tuck:13,yu:17}.

\item 
\emph{Plane and space curves} such as shape contours and motion trajectories are often compared in a group action-invariant fashion for pattern recognition. 
To compare the shape of curves with the translation and scaling effects removed, an observed smooth parametrized curve in $\bbR^d$, $d=2$ or $3$ is centered and scaled to have unit length, obtaining a centered-and-scaled curve in $\cF = \{ f:[0,1]\rightarrow\bbR^{d} \mid f(0) = 0,\, \int_{0}^{1}|f'(s)|\,ds = 1 \}$.
Each $f \in \cF$ is represented by its square root velocity function (SRVF) \citep{josh:07} $\beta_f :[0,1]\rightarrow\bbR^{d}$,
$\beta_f(s)=|f'(s)|^{-1/2}f'(s)$, $s\in [0,1]$.
Now, the SRVF $\beta_f$ lies on the Hilbert sphere $\Sinf$ in the Hilbert space  $\bbH=L^{2}([0,1],\bbR^{d})=\{h:[0,1]\rightarrow\bbR^{d}\mid\int|h(s)|^{2}\,ds<\infty\}$, equipped with the inner product $\inner{h_{1}}{h_{2}}=\int_{0}^{1}h_{1}(s)^{T}h_{2}(s)\,ds$, $h_{1},h_{2}\in\bbH$.
The spherical geometry on $\Sinf$ for the SRVFs induces a special case of the elastic metric \citep{youn:98} on the space of centered-and-scaled curves $\cF$, so that the distance between curves are given by the square root of the minimal energy to transform between them.
This geometry has demonstrated attractive and interpretable performance in practical curve matching tasks \citep{su:14,baue:17,xie:17,stra:19}. 

\end{itemize}

Applications of Hilbert sphere in computer vision, medical imaging, and biological processes necessitates hypothesis tests backed by solid theory, but there are so far no available asymptotic results to support hypothesis tests under an intrinsic geometry on $\Sinf$; 
see for example \cite{wu:14,henn:16} who applied U-statistics to perform two-sample comparisons.
For random densities, $F$-tests has been consider by \cite{pete:19-1} under a Wasserstein geometry and \cite{dube:19-1} under a more general object-oriented framework, where the latter relies on an entropy number condition that has not been verified on $\Sinf$.
This provides motivation for our work, which aims to derive asymptotic distributional results on $\Sinf$, devise valid statistical inference, and showcase applications where we performed hypothesis test for random densities.

Let $\cM$ be a metric space and $\rho: \cM \times \cM \rightarrow \bbR$ the associated distance on $\cM$.
Consider an $\cM$-valued random element $X:\Omega\rightarrow\cM$ measurable between a complete probability space $(\Omega,\cA,P)$ and $\sigma$-algebra $\mathscr{B}(\cM)$ induced on $\cM$ by $\rho$. 
The \emph{Fr\'{e}chet mean} \citep{frec:48} is a commonly used location descriptor for non-Euclidean data objects and is termed the \emph{intrinsic mean} if $\cM$ is a Riemannian manifold since it only rely on the intrinsic geometry rather than an ambient space. 
If there exists a unique minimizer of the Fr\'{e}chet functional
\begin{equation}
M(\cdot)=E\rho^{2}(X,\cdot),\label{eq:muDef}
\end{equation}
then the minimizer $\mu = \argmin_{p \in \cM} M(p)$ is called the Fr\'{e}chet mean of $X$. 
Similarly, for independent realizations $X_{1},\,\dots,\,X_{n}$ of $X$, the sample Fr\'{e}chet functional is 
\begin{equation}
M_{n}(\cdot)=\oneovern\sumin\rho^{2}(X_{i},\cdot)  \label{eq:hmuDef} 
\end{equation}
and the sample Fr\'{e}chet mean is $\hmu = \argmin_{p \in \cM} M_n(p)$ given existence and uniqueness.
This work focuses on $\cM=\Sinf$ with $\rho$ being the geodesic distance.

When the Riemannian manifold $\cM$ is finite-dimensional, theory and methods for the intrinsic Fr\'{e}chet mean have been well investigated.
Various authors have considered, for example, confidence intervals \citep{bhat:05}, hypothesis tests \citep{huck:12}, regression \citep{zhu:09}, and principal component analysis \citep{flet:04} for manifold-valued data. 
Extensive study has also characterized the existence and uniqueness \citep{karc:77,ziez:77,le:01,afsa:11}, consistency \citep{bhat:03}, and central limit theorems (CLTs) \citep{bhat:05,bhat:17} of the Fr\'{e}chet mean on a general finite-dimensional $\cM$.
A slower-than-$n^{1/2}$ CLT has been studied on the circle \citep{hotz:15} and, more recently, high-dimensional spheres and general Riemannian manifolds \citep{eltz:19}.

Also well-studied are random objects lying in an infinite-dimensional Hilbert space, which are termed functional data \citep{wang:16}.
Thanks to the flat geometry and vector space structure, the definition of the mean element is straightforward, and the asymptotic theory is obtained through extending the multivariate case \citep{hsin:15}. 
The Hilbert mean extends to a class of Hilbert manifolds through the extrinsic \citep{elli:13} or the transformation method \citep{pete:16} by mapping the original data into a Hilbert space and perform analysis there.
However, these approaches depend on the embedding or transformation and does not preserve the intrinsic geometry on $\Sinf$, so the theory derived for these methods cannot be applied to obtain an intrinsic analysis for objects lying on $\Sinf$.

Little is known about the theory for the Fr\'{e}chet mean on curved infinite-dimensional geometries.
In particular, on the Hilbert sphere $\Sinf$, the existence and uniqueness of the Fréchet mean has not been established, and no distributional results  are available for the inference of the intrinsic mean.
Difficulties in deriving these results on $\Sinf$ include the lack of compactness and the positive curvature, which prevents proof techniques developed for a finite-dimensional manifold \citep{afsa:11,bhat:05} or Hilbert space \citep{hsin:15} to be applicable. 
General results for the convergence rate of the sample Fréchet mean in abstract settings has been established by \cite{goui:19} and \cite{ahid:19,scho:19}, where the latter two works applied empirical process theory under entropy conditions, which are challenging to verify for the infinite-dimensional positively curved manifold $\Sinf$; 
no distributional results was established were provided there, and the uniqueness and existence were assumed. 

Our contributions include establishing the uniqueness and existence of the Fréchet mean on $\Sinf$ and its large sample properties, and applying these results to derive valid hypothesis tests for data objects.
We state in \ref{s:theory} theoretical properties of the intrinsic data analysis on $\Sinf$.
The existence and uniqueness of the intrinsic mean are shown in \ref{thm:muExist}, and large sample properties of the sample intrinsic mean $\hmu$ including a strong LLN and a CLT are stated in \ref{prop:muConsistent} and \ref{thm:muCLT}, respectively. 
Theoretical hurdles in verifying tightness and Lipschitz continuity are overcome with utilizing weak compactness, Hilbert differential geometry \citep{lang:99}, and a careful analysis of the spherical geometry.
The asymptotic covariance of the limiting Gaussian element in the CLT of $\hmu$ is shown to be consistently estimated by an empirical version in \ref{prop:T1A}.
To quantify the deviation of data observations around the intrinsic mean, we formulate intrinsic CLTs for the tangent vectors in \ref{cor:vCLTint} and the covariance operator in \ref{thm:covCLT} based on parallel transport.
Unlike the case in a flat Hilbert space, additional terms manifest in the asymptotic covariance of the sample covariance operator due to the curvature of $\Sinf$.
Results on the covariance will be useful to derive theoretical properties of principal component analysis \citep[e.g.][]{tuck:13} on $\Sinf$, though the latter exposition is out of scope of this work.

Asymptotic and bootstrap hypothesis tests for the population means are proposed in \ref{s:testing} as a result of the CLTs, with consistent properties of the test derived in \ref{cor:Atest1}--\ref{cor:mu2boot}. 
Two test statistics are constructed using the norm and the projections of the intrinsic sample means expressed on a chart, respectively, analogous to those developed for functional data \citep[][]{horv:13-1,aue:18}.
Simulation studies in \ref{s:simulation} demonstrate that the hypothesis tests based on the intrinsic mean have smaller bias than that based on the extrinsic mean when data are asymmetrically distributed around the mean. 
The asymptotics is shown to kick in with a moderate sample size of 50 despite the infinite-dimensional $\Sinf$. 
A real data study of daily taxi demands in Manhattan is presented in \ref{s:taxi}, where the taxi demands are modeled as densities.
The spherical geometry demonstrates advantage over alternative geometries in detecting changes in the demand patterns. 
The proofs of the main results are deferred to the Supplemental Materials.

\section{Hilbert Sphere Geometry\label{s:geometry}}

The infinite-dimensional Hilbert sphere $\Sinf$ is the unit sphere in Hilbert space $\bbH$ with well-known geometry and explicit expressions for its geometric quantities. 
In particular, the \emph{tangent space} of $\Sinf$ at $p\in\Sinf$ is the subspace $T_{p}\Sinf=\{v\in\bbH\mid\inner vp=0\}$ of $\bbH$ with codimension one. 
The \emph{metric tensor} $\innerp uv$ for tangent vectors $u$, $v\in T_{p}\Sinf$ at $p\in\Sinf$ is induced from and equal to the $\bbH$-inner product $\inner uv$; we suppress the subscript $p$ in the metric to lighten notation. %
The \emph{geodesic distance} $\rho:\Sinf\times\Sinf\rightarrow\bbR$ is given by $\rho(f,g)=\arccos(\inner fg)$. 
The Riemannian \emph{exponential map} at $p\in\Sinf$ is $\exp_{p}:T_{p}\Sinf\rightarrow\Sinf$, $\exp_{p}v=\cos(\norm v)p+\sin(\norm v)\norm v^{-1}v$, which preserves the distance to the origin, i.e. $\norm v=\rho(p,\exp_{p}v)$.
The inverse exponential map, or the \emph{logarithm map}, at $p\in\Sinf$ is defined as $\log_{p}:\Sinf\backslash\{-p\}\rightarrow T_{p}\Sinf$, $\log_{p}x=\arccos(\inner px)\norm u^{-1}u$, where $-p$ is the antipodal of $p$ and $u=x-\inner pxp$. 
A chart $\tau: U\subset \Sinf \rightarrow \bbE$ is a homeomorphism that maps $U$ onto open subset $\tau(U)$ in a Hilbert space $\bbE$.
We require $\tau$ to be smooth for computation, so that $\tau^{-1}: \tau(U) \subset \bbE \rightarrow \bbH$ is a smooth diffeomorphism. 
For example, $\tau_p(\cdot)=\log_{p}(\cdot)$ is a smooth chart defined on $U_p=\Sinf\backslash\{-p\}$, for $\pinSInf$. 
For general Hilbert differential geometry we refer to \citet{lang:99}.

For concreteness, here we consider in our numerical illustrations the Hilbert space $\bbH=L^{2}(\cI) = \{f:\cI\rightarrow\bbR\mid\int_{\cI}f(s)^{2}ds<\infty\}$ of (equivalent classes of) square-integrable functions on a compact Euclidean set $\cI$, while we note that the definition of a Hilbert sphere is fully intrinsic through isometry \citep[][]{lang:99}.

\section{Properties of the Fréchet mean on $\protect\Sinf$\label{s:theory}}

\subsection{Intrinsic mean\label{ss:mean}}

Unlike the case in a Euclidean space, the Fréchet mean for data lying on a manifold may not exist and may not be unique.  
A neighborhood condition \ref{a:muExtNeighborhood} is needed to guarantee the existence and uniqueness of the intrinsic means $\mu$ defined in~\eqref{eq:muDef}.  %
\begin{enumerate}[label=(A\arabic*)] 
\item \label{a:muExtNeighborhood} The support $\cU\subset\Sinf$ of $X$ satisfies  $P(X\in\cU)=1$ and $\sup_{x,y\in\cU}\rho(x,y)\le\pi/2$. 
\end{enumerate} 
A well-established neighborhood condition \citep{afsa:11} for a finite dimensional sphere $\cS^{d}$, $d<\infty$ is that $P(\rho(q,X)<r)=1$ for some $q\in\cS^{d}$ and $r<\pi/2$. 
A lack of compactness on the infinite-dimensional $\Sinf$ prevents the arguments of \cite{afsa:11} to apply, which utilizes compactness to show both the existence and uniqueness, the latter relying on the Poincaré--Hopf theorem for compact manifolds.
We interpret the additional concentration in \ref{a:muExtNeighborhood} as a compensation for a lack of compactness. %
The next theorem establishes the uniqueness and existence of the Fréchet mean.
\begin{thm} \label{thm:muExist} 
If \ref{a:muExtNeighborhood} holds, then there exists a unique intrinsic mean $\mu$ of $X$ on $\Sinf$. 
Furthermore, $P(\rho(\mu,X)<\pi/2)=1$. 
\end{thm} 
The existence proof is based on the weak compactness of the unit ball in $\bbH$ by the Banach--Alaoglu theorem \citep{rudi:73} and an analysis of $\rho^{2}(x,\cdot)$, 
the uniqueness result makes use of the convexity of $\rho$, 
and the proximity of $\mu$ and $X$ is obtained through a reflection argument used in \citet{afsa:11}.
With the well-definedness of the population and sample intrinsic means $\mu$ and $\hmu$ guaranteed by \ref{thm:muExist}, we derive the consistency for $\hmu$ based on $M$-estimation and empirical process \citep[see, e.g.,][]{van:96}.
\begin{prop} \label{prop:muConsistent} 
If \ref{a:muExtNeighborhood} holds, then 
\begin{equation*} 
\rho(\hmu,\mu)=o(1)\quadas\label{eq:muConsistent} 
\end{equation*} 
\end{prop} 

Some notations are needed to state the CLT for $\hmu$.
Let $\cB(E_1,E_2)$ be the Banach space of bounded linear operators between Banach spaces $E_1$ and $E_2$. 
An element $f$ in the Hilbert space $(\bbE, \innerE{\cdot}{\cdot})$ is identified with $f^{*}\coloneqq\inner f{\cdot}$ in the dual space $\bbE^{*} \coloneqq \cB(\bbE,\bbR)$ by the Riesz representation theorem.
Denote $\cA^{*}$ as the adjoint of a linear operator $\cA$ between Hilbert spaces.
The tensor product $\otimes: \bbE^* \times \bbE^* \rightarrow \cB(\bbE^{*}, \bbE^{*})$ in the dual space $E^{*}$ is defined as $(f^{*}\otimes g^{*})h^{*} = \innerE fh \, g^{*}$ for $f^{*},g^{*},h^{*}\in\bbE^{*}$. 
In a neighborhood $U$ of $\mu$, we express $\mu$, $\hmu$, and $\rho$ on a smooth chart $\tau:U\subset\Sinf\rightarrow\bbE$ and write $\mutau=\tau(\mu)$, $\hmutau = \tau(\hmu)$, and  $\rhotau:\Sinf\times\bbE\rightarrow\bbR$, $\rhotau(x,e)=\rho(x,\tau^{-1}(e))$.
Let $D_{2}$ denote the partial (Fr\'{e}chet) derivative of a multivariate function w.r.t. the second argument, where the definition of Fr\'{e}chet derivatives is reviewed in \ref{app:fret}. 
Set $\psitau:\Sinf\times\bbE\rightarrow\cB(\bbE,\bbR),$ $\psitau(x,e)=D_{2}\rhotau^{2}(x,e)$;  
$\cF_{\tau}=E[\psitau(X,\mutau)\otimes\psitau(X,\mutau)]\in\cB(\bbE^{*},\bbE^{*})$; and $\Lambda_{\tau}=ED_{2}^{2}\rhotau^{2}(X,\mutau) \in \cB(\bbE, \cB(\bbE,\bbR))$, for which the explicit forms are obtained in the Supplemental Materials. 

\begin{thm}%
\label{thm:muCLT} Let $(U,\tau)$ be a chart of $\Sinf$ in a neighborhood of $\mu$. 
Under \ref{a:muExtNeighborhood}, 
\begin{equation*} 
\sqrt{n}(\hmutau - \mutau)\conv LZ\label{eq:muCLT} 
\end{equation*} 
as $\ntoinf$, where $Z$ is a Gaussian random element in $\bbE$ with mean zero and covariance operator $\sT:\bbE\rightarrow\bbE$, $\sT=(\Lambda_{\tau}^{-1})\cF_{\tau}(\Lambda_{\tau}^{-1})^{*}$, satisfying \begin{equation*} 
\innerE{h_{1}}{\sT h_{2}}=\cov(\innerE{\Lambda_{\tau}^{-1}\psitau(X,\mutau)}{h_{1}},\innerE{\Lambda_{\tau}^{-1}\psitau(X,\mutau)}{h_{2}})\label{eq:covZ} 
\end{equation*} 
for $h_{1},\,h_{2}\in\bbE$. 
The operator $\Lambda_{\tau}$ is continuously invertible. 
\end{thm} 
\ref{thm:muCLT} follows from a linearization argument applied on the chart representation $M_{n,\tau}:\bbE\rightarrow\bbR$, $M_{n,\tau}(e)=M_{n}(\tau^{-1}(e))$ of $M_{n}(\cdot)$ in a neighborhood of $\mutau$, where the uniform convergence of the residual is obtained due to the simple dependency of the geodesic distance $\rho$ on the inner product of its arguments. 
The key difficulty is verifying the Lipschitz continuity of the criterion function under the infinite-dimensional setup by analyzing the Hilbert sphere geometry.
The CLT is intrinsic in the sense that if $\eta:V\rightarrow \bbE$ is another chart with $\mu \in V$, then $\sqrtn[\eta(\hmu)-\eta(\mu)]$ converges in law to $D(\eta\circ\tau^{-1})(\mutau)Z$ where $Z$ is the limiting Gaussian element under chart $\tau$.

The asymptotic distribution in \ref{thm:muCLT} needs to be estimated in practice for statistical inference such as in hypothesis tests that will be detailed in \ref{s:testing}. 
Define $\hsT=(\hLambda_{\tau}^{-1})\hcF_{\tau}(\hLambda_{\tau}^{-1})^{*}$, $\hcF_{\tau}=n^{-1}\sum_{i=1}^{n}\psitau(X_{i},\hmu_{\tau})\otimes\psitau(X_{i},\hmu_{\tau})$, and $\hLambda_{\tau}=n^{-1}\sum_{i=1}^{n}D_{2}^{2}\rhotau^{2}(X_{i},\hmu_{\tau})$. 
For an operator $\cC:\bbF\rightarrow\bbF$ defined on a Hilbert space $\bbF$, let $\norm{\cC}=\sup_{\norm h_{\bbF}=1}\norm{\cC h}_{\bbF}$ denote the operator norm and $\norm{\cC}_{1}=\sumjinf\inner{e_{j}}{(\cC^*\cC)^{1/2} e_{j}}_\bbF$ the trace norm, where $\{e_j\}_{j=1}^\infty$ is an arbitrary complete orthonormal basis of $\bbF$.
The next result states that the asymptotic distribution of the limiting
Gaussian element can be estimated consistently.
\begin{prop} %
\label{prop:T1A} Under the conditions of \ref{thm:muCLT},
as $n\toinf,$
\[
\norm{\hsT-\sT}_{1}=o(1) \quadas
\]
Let $Z_{n}$ and $Z$ be zero-mean Gaussian random elements in $\bbE$ with covariance $\hsT$ and $\sT$, respectively. 
Then 
\[
Z_{n}\conv LZ
\]
as $\ntoinf$.
\end{prop}

\subsection{Tangent Vectors and the Covariance Operator\label{ss:tangentCov}}

On a nonlinear manifold, the deviations of data observations around the intrinsic mean are commonly characterized by the corresponding tangent vectors constructed from the logarithm maps.
The tangent vector of $X$ and its empirical version are respectively defined as 
\begin{equation*} 
V=\logmu X,\quad\hV=\loghmu X.\label{eq:VDef} 
\end{equation*} 
Similarly, let $V_{i}= \logmu X_{i}$ and $\hV_{i} = \loghmu X_{i}$ denote the corresponding quantities for an observation $X_{i}$, $i=1,\dots,n$. 
\begin{prop} \label{prop:Vmean} 
If $\mu$ is the intrinsic mean of an $\Sinf$-valued random variable $X$ with $P(\rho(\mu,X)=\pi)=0$, then $EV=0$. 
\end{prop} 
\ref{prop:Vmean} states that the logarithm map centers the observations around the intrinsic mean on $\Sinf$, a property that has been established on finite-dimensional Riemannian manifolds \citep{karc:77,bhat:03}.

The tangent space $\TmuS$ at $\mu$ is a Hilbert space with inner product given by the metric tensor.
We obtain first a CLT of $\hmu$ on the tangent space $\TmuS$ as a corollary of \ref{thm:muCLT} by considering $\tau=\log_{\mu}(\cdot)$, so that $\psitau(X, \mutau) = -2\inner{V}{\cdot}$ as will be shown in the proof of \ref{prop:Vmean}. %
Here $\Lambda_\tau$ becomes $\Lambda=ED_{2}^{2}\rho^{2}(X,\mu) \in \cB(T_{\mu}\Sinf, \cB(T_{\mu}\Sinf,\bbR))$, which is identified with the $T_{\mu}\Sinf$-valued linear map $\Lambda_{1} \in \cB(T_{\mu}\Sinf, T_{\mu}\Sinf)$ so that $\Lambda v = \inner{\Lambda_{1}v}{\cdot}$, for $v\in T_{\mu}\Sinf$.
\begin{cor}
\label{cor:muCLTlog} Under the conditions of \ref{thm:muCLT}, 
$\sqrtn\log_{\mu}\hmu$ converges in distribution to a Gaussian random element on $T_{\mu}\Sinf$ with mean zero and covariance $4\Lambda_{1}^{-1}E[V \otimes V ]\Lambda_{1}^{-1}$.
\end{cor}

The variation of $X$ around the intrinsic mean is summarized by the covariance operator, which is a linear characterization of the intrinsic variation and has been applied to generalize the principal component analysis to Riemannian manifolds \citep{flet:04,laza:17}.
The population and empirical covariances are $\cG: \TmuS\rightarrow \TmuS$, $\cG=E(V\otimes V)$ and $\hcG:\ThmuS \rightarrow \ThmuS$, $\hcG=n^{-1} \sumin\hV_{i}\otimes\hV_{i}$, respectively, where $\otimes$ is the tensor product on the tangent spaces, so that 
\begin{equation*}
\cG h=E(\inner VhV), \quad\, \hcG g=\oneovern\sumin \inner{\hV_{i}}g\hV_{i} \label{eq:hcGDef}
\end{equation*}
for $h\in\TmuS$ and $g\in\ThmuS$. 

A corollary for comparing tangent vectors is needed for assessing the estimation of the covariance operator.
Parallel transport under the Levi-Civita connection on $\Sinf$ \citep{lang:99} is applied to perform intrinsic comparisons of tangent vectors on different tangent spaces. 
Let $P_{x}^{y}:T_{x}\Sinf\rightarrow T_{y}\Sinf$ denote the parallel transport of tangent vectors on $T_x\Sinf$ along the geodesic leaving from $x$ to $y$ on $\Sinf$, where $x$ and $y$ are not antipodal so that the geodesic is unique. 
Write  $\Lambda_{x}=D_{2}^{2}\rho^{2}(x,\mu)\in\cB(T_{\mu}\Sinf,\cB(T_{\mu}\Sinf,\bbR))$ and we identify it with $\Lambda_{x1}\in\cB(T_{\mu}\Sinf,T_{\mu}\Sinf)$ such that $\Lambda_{x}v = \inner{\Lambda_{x1}v}{\cdot}$ for any $v\in T_{\mu}\Sinf$. 
Recall that $\Lambda_{1}$ is defined before \ref{cor:muCLTlog}.

\begin{cor} %
\label{cor:vCLTint} Under \ref{a:muExtNeighborhood}, $v=\logmu x$
and $\hv=\loghmu x$ are well defined for any $x\in\Sinf$ with $\rho(x,\mu)<\pi$, the latter almost surely as $\ntoinf$. 
Moreover, 
\begin{equation*}
\sqrt{n}(P_{\hmu}^{\mu}\hv-v)\conv LZ_{v}
\end{equation*}
as $\ntoinf$, where $Z_{v}$ is a Gaussian random element in $T_{\mu}\Sinf$ with
mean 0 and covariance $\Lambda_{x1}\Lambda_{1}^{-1}E[V\otimes V]\Lambda_{1}^{-1}\Lambda_{x1}$.
\end{cor}

To derive a central limit theorem for the covariance, operators $\cG$ and $\hcG$ are analyzed as Hilbert--Schmidt operators on the tangent spaces and are compared through the parallel transport.
For $x\in\Sinf$ lying in a small neighborhood of $\mu$, let $\cP_{x}^{\mu}:\cB(T_{x}\Sinf,T_{x}\Sinf)\rightarrow\cB(T_{\mu}\Sinf,T_{\mu}\Sinf)$ be the parallel transport of operators from $x$ to $\mu$ such that for any operator $\cA\in\cB(T_{x}\Sinf,T_{x}\Sinf)$, $\cP_{x}^{\mu}(\cA)v=P_{x}^{\mu}\cA P_{\mu}^{x}v$, $v\in T_{\mu}\Sinf$. 
For a Hilbert space $\bbF$, let $\cB_{\HS}(\bbF,\bbF)$ denote the Hilbert space of Hilbert--Schmidt operators on $\bbF$ equipped with the inner product %
\[
\innerHS{\cF_{1}}{\cF_{2}}=\sumjinf\inner{\cF_{1}e_{j}}{\cF_{2}e_{j}}_\bbF,
\]
where $\cF_{1}$, $\cF_{2}$ are operators on $\bbF$ with finite Hilbert--Schmidt norm induced by this inner product, and $\{e_{j}\}_{j=1}^{\infty}$ is a complete orthonormal basis of $\bbF$. 
Also let $\otimes_{\HS}$ denote the tensor product in $\cB_{\HS}(\TmuS,\TmuS)$.
Define $\cH:\TmuS \rightarrow\cB_{\HS}(\TmuS,\TmuS)$, $\cH(\cdot)=E[\Lambda_{X1}(\cdot) \otimes V]$ where $\Lambda_{X1}$ is $\Lambda_{x1}$ at a random $x=X$. %

\begin{thm} %
\label{thm:covCLT}
Under \ref{a:muExtNeighborhood}, 
\begin{equation*}
\sqrt{n}(\cP_{\hmu}^{\mu}\hcG-\cG)\conv LZ_{\cG}\label{eq:covCLT}
\end{equation*}
in $\cB_{HS}(T_{\mu}\Sinf,T_{\mu}\Sinf)$, where $Z_{\cG}$ is a Gaussian
random element with mean 0 and covariance $E[ \cG_0 \otimes_{\HS} \cG_0 ]$, where 
\begin{equation*}
\cG_0 = V\otimes V-\cG+\cH\Lambda_1^{-1}V+(\cH\Lambda_1^{-1}V){}^{*}.  \label{eq:covCLTG0}
\end{equation*}
\end{thm}
\ref{thm:covCLT} is an extension of the central limit theorem for the covariance operator in $\bbH$ \citep[e.g. Theorem 8.1.2 in][]{hsin:15} to the Hilbert sphere $\Sinf$, where in the former Hilbert space the parallel transport is simply the identity. 
Additional terms involving $\cH$ appear in the random Hilbert--Schmidt operator $\cG_0$ that generates the asymptotic covariance due to the positive curvature of $\Sinf$.
\ref{thm:covCLT} can be applied to derive asymptotic theory for the principal component analysis \citep{tuck:13,dai:18-5,lin:19} based on $\hcG$, an exposition that is beyond the scope of this work.

\section{Hypothesis Tests for the Intrinsic Mean}
\label{s:testing}

\subsection{General Setup}
\label{ss:testSetup}

We obtain here one- and two-sample hypothesis tests of the intrinsic mean as a result of the CLT in \ref{thm:muCLT}; for completeness, a paired-sample test is formulated in the Supplemental materials. 
To handle the manifold constraint, these tests make use of a chart $(U, \tau)$ so that test statistics based on the intrinsic mean can be constructed analogously to those in a Hilbert space \citep{berk:09,aue:18}. 

\subsubsection*{One-sample test}

Given a sample $X_{1},\,\dots,\,X_{n}$ on $\Sinf$ with unknown mean $\mu$, consider the one-sample hypothesis
\begin{gather*}
H_{0}:\mu=\mu_{0},\quad H_{1}:\mu\ne\mu_{0},\label{eq:h1}
\end{gather*}
for a pre-specified $\mu_0 \in \Sinf$.
We propose a norm-based and a projection-based test statistic defined as, respectively, 
\begin{gather*}
T_{1}=n\normE{\tau(\hmu)-\tau(\mu_{0})}^{2}, \label{eq:T1}\\
S_{1}=n\sumkK\frac{\innerE{\tau(\hmu)-\tau(\mu_{0})}{\hphi_{k}}^{2}}{\hlambda_{k}},\label{eq:S1}
\end{gather*}
where the projection-based test utilizes $K < \infty$ projections.
Here $S_1$ is the empirical estimate of $\tilde{S}_{1}=  n\sumkK\lambda_{k}^{-1}\innerE{\tau(\hmu)-\tau(\mu_{0})}{\phi_{k}}^{2}$, 
and the $(\lambda_{k},\phi_{k})$ and $(\hlambda_k, \hphi_k)$ are the eigenvalue--eigenfunction pairs of the true and estimated covariance operator $\sT$ or $\hsT$ of the limiting element $Z$ in \ref{thm:muCLT}, respectively. 
Note that one cannot obtain a Hotelling's $T^2$-like test statistic by normalizing the mean with the inverse of its covariance operator since the resulting test statistic diverges \citep{haje:62}.
A natural choice of chart is $\tau(\cdot)=\log_{\mu_{0}}(\cdot)$, under which the norm-based test statistic reduces to $T_{1}=n\rho^{2}(\hmu,\mu_{0})$.

By Slutsky's theorem and \ref{thm:muCLT}, the limiting distribution for $T_{1}$ is $\normE{Z}^2$ and that for $S_1$ is $\chi_{K}^{2}$. 
Summarizing,

\begin{cor}
\label{cor:Atest1}Assume the conditions of \ref{thm:muCLT} hold.
Then
\begin{align*}
T_{1} & \conv L\sumkinf\lambda_{k}W_{k},\\
S_{1} & \conv L\chi_{K}^{2}, 
\end{align*}
where $W_{k}$ are i.i.d. $\chi_{1}^{2}$ random variables.
\end{cor}

\subsubsection*{Two-sample test}
The two-sample setting is that we have i.i.d. observations $X_{1}^{[g]},\,\dots,\,X_{n_{g}}^{[g]}$ on $\Sinf$ from Population~$g$, $g=1,2$, where $n_g$ is the sample size and $n=n_1 + n_2$ is the total sample size.
The unknown intrinsic mean in Population~$g$ is $\mu_g$ and is estimated by the empirical mean $\hmu_g$. 
The two-sample hypothesis is
\begin{equation*}
H_{0}:\mu_{1}=\mu_{2},\quad H_{1}:\mu_{1}\ne\mu_{2}.\label{eq:h2}
\end{equation*}
Under chart $\tau$, let $\sT_{g}=(\Lambda_{g\tau}^{-1})\cF_{g\tau}(\Lambda_{g\tau}^{-1})^{*}$ be the asymptotic covariance of $\tau(\hmu_{g})$ given by \ref{thm:muCLT},
where $\cF_{g\tau}=E[\psitau(X_{i}^{[g]},\tau(\mu_{g}))\otimes\psitau(X_{i}^{[g]},\tau(\mu_{g}))]$,
$\Lambda_{g\tau}=ED_{2}^{2}\rhotau^{2}(X_{i}^{[g]},\tau(\hmu_{g}))$,
$g=1,\,2$. The pooled covariance operator is written as $\sT_{\pool}=(n/n_{1})\sT_{1}+(n/n_{2})\sT_{2}$.
The following corollary of \ref{thm:muCLT} provides a basis for the
two-sample tests.
\begin{cor}
\label{cor:mu2dist} Assume $n_{1},n_{2}\toinf$ and $n_{1}/n_{2}\rightarrow q$ for some $q\in(0,1)$. 
If $H_{0}$ is true and the conditions of \ref{thm:muCLT} hold for both populations, then
\[
\sqrt{n}[\tau(\hmu_{1})-\tau(\hmu_{2})]\conv LZ_{2},
\]
where $Z_{2}$ is a Gaussian random element with mean $0$ and covariance $\sT_{\pool}$.
\end{cor}
The norm-based and projection-based two-sample test statistics are
\begin{gather*}
T_{2}=n\normE{\tau(\hmu_{1})-\tau(\hmu_{2})}^{2},\label{eq:T2}\\
S_{2}=n\sumkK\frac{\innerE{\tau(\hmu_{1})-\tau(\hmu_{2})}{\hphi_{k,\pool}}^{2}}{\hlambda_{k,\pool}}. \label{eq:S2}
\end{gather*}
Here $(\hlambda_{k,\pool},\hphi_{k,\pool})$, $k=1,2,\dots$ are the eigenvalue--eigenfunction pairs of the estimated pooled covariance $\hsT_{\pool}=(n/n_{1})\hsT_{1}+(n/n_{2})\hsT_{2}$,
where $\hsT_{g}=(\hLambda_{g\tau}^{-1})\hcF_{g\tau}(\hLambda_{g\tau}^{-1})^{*}$,
$\hcF_{g\tau}=n_{g}^{-1}\sum_{i=1}^{n_{g}}\psitau(X_{i}^{[g]},\tau(\hmu_{g}))\otimes \psitau(X_{i}^{[g]},\tau(\hmu_{g}))$,
and $\hLambda_{g\tau}=n_{g}^{-1}\sum_{i=1}^{n_{g}}D_{2}^{2}\rhotau^{2}(X_{i}^{[g]},\tau(\hmu_{g}))$, $g=1,2$. 
\begin{cor}
\label{cor:Atest2}Assume the conditions of \ref{cor:mu2dist} hold.
Let $W_{k}$ be i.i.d. $\chi_{1}^{2}$ random variables and $\lambda_{k,\pool}$, $k=1,2,\dots$ be the eigenvalues of $\sT_{\pool}$. Then
\begin{align*}
T_{2} & \conv L\sumkinf\lambda_{k,\pool}W_{k},\\
S_{2} & \conv L\chi_{K}^{2}.
\end{align*}
\end{cor}

\subsection{Asymptotic Tests\label{ss:atest}}

In performing asymptotic tests, the asymptotic null distributions of the norm-based test statistics involve unknown eigenvalues, which need to be estimated. 
The infinite sum in the limiting distribution poses the question of whether the asymptotic $p$-values can be estimated consistently. 
The answer is affirmative: 
In the one-sample scenario, the asymptotic null distribution of $T_{1}$ is $\norm{Z}^2$, which is estimated by the distribution of $\hat{T}_1 \coloneqq \norm{Z_n}^2 = \sumkinf\hlambda_{k}W_{k}$.
As a consequence of \ref{prop:T1A} and the continuous mapping theorem, the limiting distributions of $T_1$ and $\hat{T}_1$ are the same, so the quantiles of $T_1$ can be estimated consistently by those of $\hat{T}_1$ in the large sample limit. 
The asymptotic $p$-values can thus be consistently estimated by the quantiles of $\hat{T}_1$, defining a valid asymptotic test. 
In practice, quantiles of $\hat{T}_1$ are obtained by Monte Carlo. 
Although the asymptotic distribution is non-pivotal, the \ref{thm:muCLT} leads to the explicit form of the asymptotic distribution of the test statistics, enabling efficient implementation. 

For the projection-based tests, the tuning parameter $K$ needs to be chosen in order to maximally capture the difference in the means. 
While an optimal choice may depend on the sample size and the stochastic structure of the data, 
we find in our numerical studies that the Fraction of Variance Explained (FVE) criterion leads to reasonable test performance,  by setting $K = K^*$, 
\begin{equation}
K^{*}=\min\{K\ge1\mid\frac{\sum_{j=1}^K \hlambda_{j}}{\sumkinf\hlambda_{k}}\ge r\}, \label{eq:FVE}
\end{equation}
with a threshold $r\in(0,1)$. 
The projection-based tests will be powerful as long as the mean difference (on the chart) is not orthogonal to the subspace spanned by the first $K$ eigenfunctions. 
Our experience, in agreement with \cite{berk:09,horv:13-1}, is that in practice the mean difference is usually well captured by the first few projections. 
Even though the norm-based tests always capture the mean difference, 
the projection-based tests are often more powerful in our numerical studies since they utilize a pivotal test statistic and focus on the leading components with higher signal-to-noise ratio. 
The asymptotic tests in the two- and paired-sample scenarios follow analogous development. 

\subsection{Bootstrap Tests\label{ss:btest}}

Alternative bootstrap tests are proposed to improve finite sample performance. 
The bootstrap tests do not require the computation of the covariance components of the norm-based test statistics, and enjoy a second order accuracy using the pivotal projection-based test statistics \citep{hall:92-2}. 

In the one-sample scenario, let $X_{1}^{*},\,\dots,\,X_{n}^{*}$ be a nonparametric bootstrap sample drawn from $X_{1},\,\dots,\,X_{n}$ with replacement and $\hmu^{*}$ be the bootstrap sample intrinsic mean.
The bootstrap distributions are derived from 
\begin{align*}
T_1^* & = n \normE{\tau(\hmu^*) - \tau(\hmu)}^2, \\
S_1^* & = n \sumkK \frac{\innerE{\tau(\hmu^*)-\tau(\hmu)}{\hphi_{k}^*}^{2}}{\hlambda_{k}^*},
\end{align*}
where the $(\hlambda_k^*, \hphi_k^*)$ are the eigenpairs of $\hsT^*_\tau$ constructed analogously to the sample covariance $\hsT_\tau$ but with the bootstrap sample.
The validity of the bootstrap test is guaranteed by the following corollary of \ref{thm:muCLT} and the bootstrap theorem for infinite-dimensional $Z$-estimators in \citet{well:96}. 
\begin{cor}%
\label{cor:mu1boot} Under the conditions of \ref{thm:muCLT}, the bootstrap intrinsic mean $\hmu^{*}$ is consistent for $\mu$. Furthermore, as $\ntoinf$, 
\[
\sqrt{n}[\tau(\hmu^{*})-\tau(\hmu)]\conv LZ^{*},
\]
where $Z^{*}$ is a Gaussian random element sharing the same distribution with $Z$ in \ref{thm:muCLT}. 
Hence, the asymptotic distributions are the same for $T_{1}^{*}$ and $T_1$, as well as for $S_{1}^{*}$ and $S_{1}$.
\end{cor}

For the two-sample test, let $X_{1}^{[g]*},\,\dots,\,X_{n_{g}}^{[g]*}$ be a nonparametric bootstrap sample of size $n_{g}$ from Population~$g$, $g=1,2$. 
Let $\hmu_{g}^{*}$ be the intrinsic mean and $\hsT_{g}^{*}$ the sample covariance operator of the bootstrap sample in Population~$g$, $\hsT_{\pool}^{*}=(n/n_{1})\hsT_{1}^{*}+(n/n_{2})\hsT_{2}^{*}$ the bootstrap pooled covariance, and $(\hlambda_{k,\pool}^{*},\hphi_{k,\pool}^{*})$, $k=1,2,\dots$ the eigenpairs of $\hsT_{\pool}^{*}$. 
The bootstrap versions of for $T_{2}$ and $S_{2}$ are, respectively, 
\begin{align*}
T_{2}^{*} & =n\normE{[\tau(\hmu_{1}^{*})-\tau(\hmu_{1})]-[\tau(\hmu_{2}^{*})-\tau(\hmu_{2})]}^{2},\\
S_{2}^{*} & =n\sumkK\frac{\innerE{[\tau(\hmu_{1}^{*})-\tau(\hmu_{1})]-[\tau(\hmu_{2}^{*})-\tau(\hmu_{2})]}{\hphi_{k,\pool}^{*}}^{2}}{\hlambda_{k,\pool}^{*}}, %
\end{align*}
where the bootstrap statistics are constructed as such to increase
power.
\begin{cor} %
 \label{cor:mu2boot}Under the conditions of \ref{cor:mu2dist},
the bootstrap intrinsic mean $\hmu_{g}^{*}$ is consistent for $\mu_{g}$, $g=1,\,2$. 
Moreover,
\[
\sqrt{n}\left\{ [\tau(\hmu_{1}^{*})-\tau(\hmu_{1})]-[\tau(\hmu_{2}^{*})-\tau(\hmu_{2})]\right\} \conv LZ_{2}^{*},
\]
where $Z_{2}^{*}$ is a Gaussian random element sharing the same distribution as $Z_{2}$ in \ref{cor:mu2dist}. 
Hence, the asymptotic distributions are the same for $T_{2}^{*}$ and $T_{2}$, as well as for $S_{2}^{*}$ and $S_{2}$.
\end{cor}

\section{Simulation Studies\label{s:simulation}}

Simulation studies were performed on $\Sinf$ modeled in $\bbH=L^{2}([0,1])$ to demonstrate the numerical properties of the proposed hypothesis tests. 
We focus on the two-sample case here, while analogous results for the one-sample scenario are included in the Supplemental Materials. 
In the two-sample scenario, we set the intrinsic population mean $\mu_{1}$ to the square root of the Beta$(2,1)$ density function and $\mu_{2}=\exp_{\mu_{1}}\delta v$, where $\delta\in[-0.4,0.4]$ is the effect size, $v=K_{\mu}^{-1/2}\sum_{k=1}^{K_{\mu}}\phi_{1k}$, $K_{\mu}\in\{1,3,5\}$ is the number of mean components, and the $\phi_{1k}$ are orthonormal functions to be detailed shortly. 
Independent observations in Population $g = 1,2$ were generated according to $X_{i}^{[g]}(s)=\exp_{\mu_{g}}\{(-1)^{g-1}\sum_{k=1}^{K_X}\xi_{gik}\phi_{gk}(s)\}$, $s\in[0,1]$ with $K_X=50$ components. 
The $k$th scores $\xi_{gik}$ were i.i.d. real-valued random variables with mean 0 and variance $\theta_k = 3^{-k}$, generated from either the normal distribution or the centered exponential distribution, i.e.  $\xi_{gik}=\eta_{gik}-E\eta_{gik}$ and $\eta_{gik}$ follows Exponential$(\theta_{k})$, for $i=1,\dots,n_g$, $g=1,2$. 
The orthonormal basis functions $\phi_{gk}$ were defined as $\phi_{gk}=R_{\mu_{g}}(\psi_{k+1})$, $k=1,\,\dots,\,K_X$, where
$\{\psi_{j}\}_{j=1}^{\infty}=$ $\{\psi_{1}(s)=1,\,\psi_{2k}(s)=2^{1/2}\sin(2k\pi s),\,\psi_{2k-1}(s)=2^{1/2}\cos(2(k-1)\pi s),\,s\in[0,1], \text{ for } k\in\bbN\}$ is the trigonometric basis on $[0,1]$, 
and $R_{q}:\bbH\rightarrow\bbH$ is the rotation operator from $\psi_{1}$ to $q \ne -\psi_1$ along the shortest geodesic, defined by %
\begin{equation*}
R_{q}(p)=p+\sin(\rho_{q})(\inner{u}p q -  \inner{q}p u) + (\cos(\rho_{q})-1)(\inner{q}p q + \inner{u}p u),\label{eq:rotation}
\end{equation*}
where $\rho_{q} = \rho(\psi_{1},q)$ and  $u=(\psi_{1}-\inner{\psi_{1}}{q}q)/(1-\inner{\psi_{1}}{q}^{2})^{1/2}$.
The mean squared geodesic distance from the observations to the intrinsic mean was $\sum_{k=1}^{K_X} \theta_{k}=0.5$, 
and the cumulative FVE by the first $J=1,\dots,5$ components, defined as $\text{FVE}(J)=\sum_{j=1}^{J}\theta_{j}/\sum_{k=1}^{K_X} \theta_{k}$, were 66.7\%, 88.9\%, 96.3\%, 98.8\%, and 99.6\%, respectively.

The simulation settings consisted of all combinations of 
\begin{itemize}
\item sample size $n_g \in \{10,\, 25,\, 50\}$;
\item number of components $K_{\mu}\in\{1,3,5\}$ in the mean difference; and 
\item either symmetric or asymmetric data generated around the mean, corresponding to the normally (norm) and exponentially (exp) distributed $\xi_{gik}$, respectively.
\end{itemize}
We compared the proposed asymptotic and bootstrap tests which are intrinsic to $\Sinf$, as well as a norm-based bootstrap test of the extrinsic means \citep{elli:13} in the ambient space $\bbH$ projected back onto $\Sinf$. 
The number of components $K$ for our projection-based tests were chosen according to the FVE criterion (\ref{eq:FVE}) with threshold $r=0.8,\,0.95$, or $0.99$, which roughly correspond to $K=2,\,3$, and 4 projections in our settings, respectively.  

\begin{sidewaysfigure}[H]
\includegraphics[width=1\textwidth]{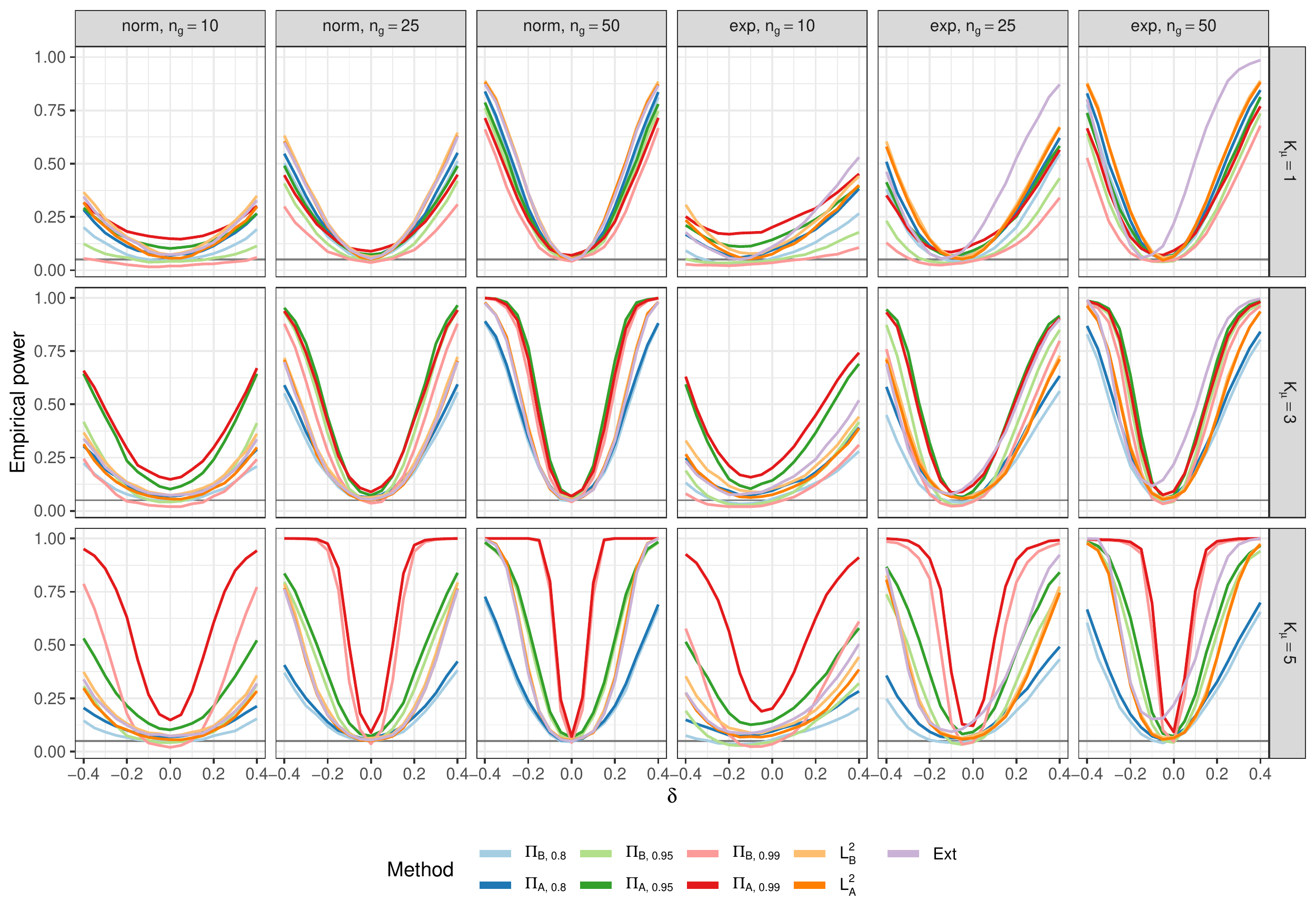}\caption{\label{fig:test2}
Empirical power curves for the two-sample tests. 
Columns correspond to different generating distributions for the principal components $\xi_{gik}$ and sample sizes $n_g$, and rows correspond to different numbers of components $K_\mu$ in the mean difference. 
The horizontal gray lines indicate the nominal level $\alpha=0.05$. 
$\Pi_{A,r}$ and $\Pi_{B,r}$, projection-based tests with FVE threshold $r$; $L_{A}^{2}$ and $L_{B}^{2}$, norm-based tests; Ext, the extrinsic bootstrap test. 
Subscripts $A$ and $B$ stand for the proposed tests in the asymptotic and bootstrap versions, respectively.}
\end{sidewaysfigure}

Calculated with 1000 Monte Carlo iterations each with 499 bootstrap samples, the empirical power curves for the two-sample tests over effect size $\delta \in [-0.4, 0.4]$ are displayed in \ref{fig:test2}, noting that $H_0$ holds if and only if $\delta = 0$.
As a visual aid, dark and light paired colors were used to denote the proposed asymptotic and bootstrap tests, respectively. 
Bootstrap tests were overall more conservative and had better control of the type I error rate (size) than the asymptotic tests, which is most apparent for $n_{g}=10$ or $25$. 
When $n_{g}=10$, the asymptotic projection-based tests were over-liberal, while the corresponding bootstrap tests properly controlled the size to be around the nominal level $\alpha=0.05$.
This can be attributed to the second-order correctness for the bootstrap tests when the test statistic is pivotal \citep{hall:92-2}, of which the effect is prominent in small samples. 
The asymptotic and bootstrap tests had almost identical performance under $n_g= 50$ (3rd and 6th columns, \ref{fig:test2}), showing that the asymptotics comes into force under this moderate sample size even if the data lie on an infinite-dimensional curved manifold $\Sinf$. 

When $n_{g}=10$, the norm-based tests had higher power than the projection-based tests that respect the nominal level.
In this small-sample scenario, the norm-based tests avoid estimating the projection directions and gain stability as compared to the projection-based tests.
With $n_g = 25$ or $50$, the projection-based tests were more powerful than the norm-based test when larger FVE thresholds were chosen to capture mean differences in multiple projections. 
Specifically, when $K_{\mu}=3$ and 5, the projection-based tests with $\text{FVE}=0.95$ and $0.99$ were the most powerful, respectively. 
This is due to the fact that the projection-based tests focus on the mean differences only in the directions with high signal-to-noise ratios. 
When $K_\mu = 1$, the norm-based tests and the projection-based tests with $\text{FVE}=0.80$ had similar performance and were both among the best performers.

The extrinsic and intrinsic norm-based tests performed similarly in the symmetric scenarios (1st--3rd columns, \ref{fig:test2}) but rather differently in the asymmetric scenarios (4th--6th columns, \ref{fig:test2}). 
All tests suffered from finite-sample biases to different extents in the asymmetric scenarios, which is reflected by the lack of power when $\delta$ was slightly below $0$.
The bias for the proposed intrinsic methods reduced as $n_{g}$ increased in the asymmetric scenarios, reaching near-unbiasedness when $n_g = 50$, but that for the extrinsic tests remained even with $n_{g}=50$.
This underlines the importance of respecting the intrinsic geometry when data is asymmetrically generated on the manifold.

\section{Data Application: Taxi Demands\label{s:taxi}}

Better understanding of the taxi demands will provide key insight into a more reliable and economic public transportation infrastructure. 
This subject has been of increasing modeling interest \citep{chu:19,dube:20} as the popularity of app-based for-hire vehicle services such as Uber and Lyft rises and data become available. 
We analyzed the demand patterns of taxi and other for-hire vehicles in New York City, 
which were extracted from a total of 1.1 billion trips records of for-hire vehicles including the yellow and green cabs, Uber, Lyft, etc. 
The trip record data were made public following Freedom of Information Law (FOIL) requests, and our analysis was built upon a database compiled by Todd Schneider available on \url{https://github.com/toddwschneider/nyc-taxi-data}.

\begin{figure}[H]
\centering
\includegraphics[width=.8\textwidth]{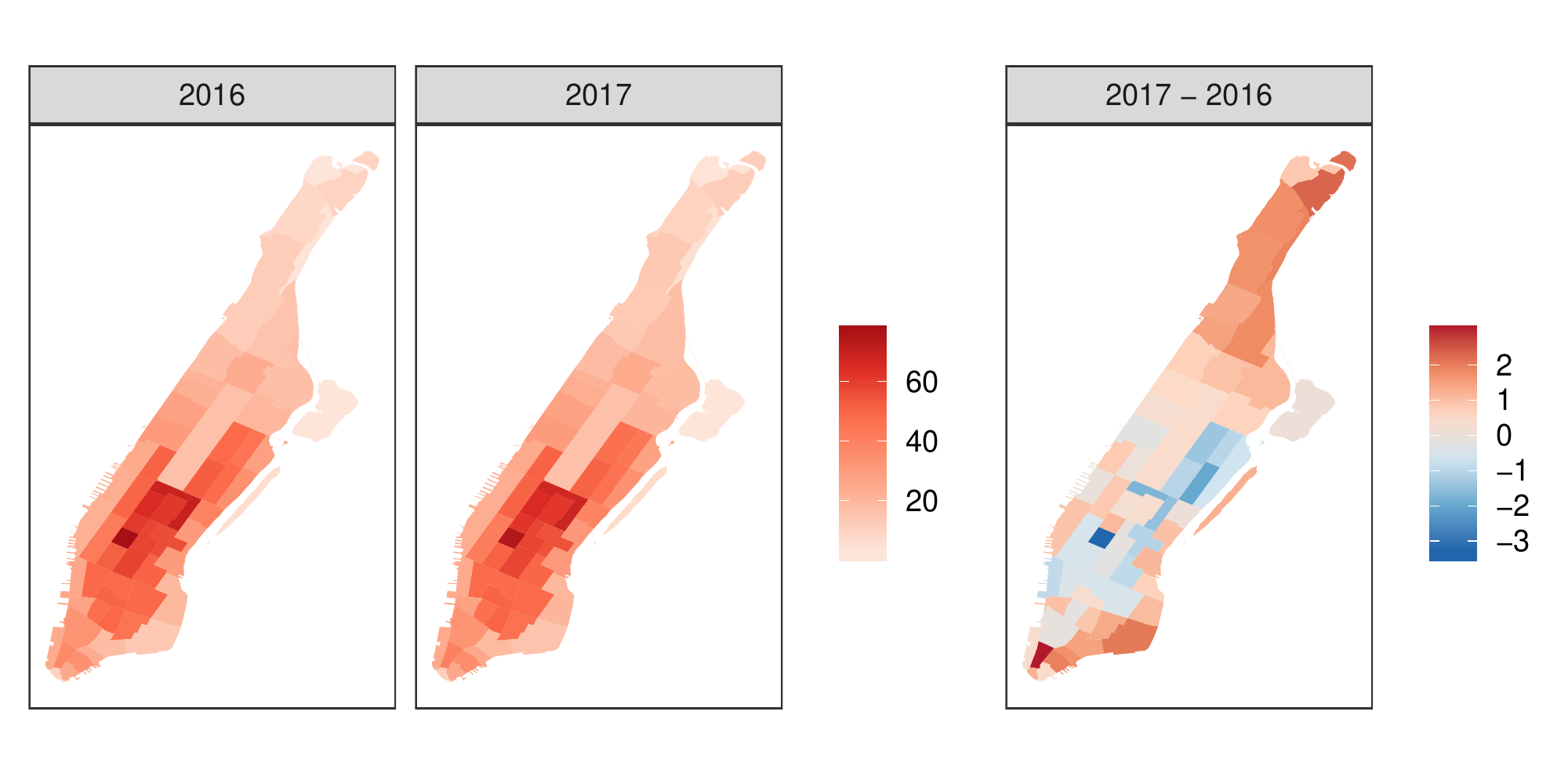}\caption{\label{fig:taxi}The sample intrinsic means and their (Euclidean) difference.}
\end{figure}

The daily taxi demand is modeled as a spatial density of the passenger pick-up locations.
For the interest of monitoring evolving demand patterns, our goal here is to compare the taxi demands in the year of 2016 and 2017. 
We focused on 
the Manhattan taxi zones $\cS$ and measured the daily demand by the probability density function $Y_{gi}(s)$ of pick-up locations in the $i$th day of year $g\in \{2016,2017\}$, where $s\in\cS$ stands for the taxi zone. 
Each spatial demand density is then transformed into a square root density $X_{gi}$ according to $X_{gi}(s)=\sqrt{Y_{gi}(s)}$, $s\in\Sinf$, so that $\norm{X_{gi}}^2 = \int_{\cS} X_{gi}^2(s) ds = \int_{\cS}Y_{gi}(s)ds = 1$ and thus $X_{gi}$ lies on the unit Hilbert sphere $\Sinf$. 
Though the average demands in $2016$ and $2017$ as measured by the intrinsic means of the square root densities $X_{gi}$ within each year appear overall similar (left two panels, \ref{fig:taxi}), there is a substantial decrease in demand concentrated near the Upper East Side of Manhattan (right panel, \ref{fig:taxi}), which is probably due to the opening of three Second Avenue subway stations on January 1, 2017. %

We investigated the effect of geometry for detecting changes in taxi demands by comparing the size and power of two-sample hypothesis tests.
The proposed intrinsic tests and the extrinsic test were applied on the square root densities $X_{gi}$ utilizing the spherical geometry, and an alternative norm-based bootstrap test was performed for the original densities $Y_{gi}$ as elements in $L^{2}(\cS)$, assuming a flat geometry. 
Two scenarios were considered, namely an Equal Mean scenario where the observations for both populations were randomly sampled without replacement from year 2017, representing that $H_0$ holds true, and an Unequal Mean scenario where the observations for the two populations were sampled from 2016 and 2017, respectively, representing $H_1$. 
The number of samples $n_{g}$ from each population varied among $10,\,15,\,20$ and $25$. 
The number of projections for the projection-based tests were  selected according to the FVE criterion with threshold $r=0.8,\,0.9$, and 0.95, which corresponds to $K=3$, 5, and 9 components, respectively. 
The empirical sizes and powers are reported in \ref{tab:taxi} for the nominal level $\alpha=0.05$, calculated from 2000 Monte Carlo iterations and 999 bootstrap samples.

\begin{table}[H]
\caption{\label{tab:taxi}
Proportions of rejected $H_{0}$ for the taxi demands at the nominal level $\alpha=0.05$. 
$\Pi_{A,r}$ and $\Pi_{B,r}$, projection-based tests with FVE threshold $r$; $L_{A}^{2}$ and $L_{B}^{2}$, norm-based tests; Ext, the extrinsic bootstrap test on $\protect\Sinf$; Dens, the bootstrap test using the original densities. Subscripts $A$ and $B$ stand for the proposed asymptotic and bootstrap tests, respectively.}

\setlength{\tabcolsep}{5pt}

\begin{tabular}{c|r@{\extracolsep{0pt}.}lr@{\extracolsep{0pt}.}lr@{\extracolsep{0pt}.}lr@{\extracolsep{0pt}.}l|r@{\extracolsep{0pt}.}lr@{\extracolsep{0pt}.}lr@{\extracolsep{0pt}.}lr@{\extracolsep{0pt}.}l|r@{\extracolsep{0pt}.}l|r@{\extracolsep{0pt}.}l}
\multicolumn{21}{l}{\emph{  Equal mean ($H_0$)}}\tabularnewline
$n_{g}$ & \multicolumn{2}{c}{$L_{A}^{2}$} & \multicolumn{2}{c}{$\Pi_{A,0.8}$} & \multicolumn{2}{c}{$\Pi_{A,0.9}$} & \multicolumn{2}{c|}{$\Pi_{A,0.95}$} & \multicolumn{2}{c}{$L_{B}^{2}$} & \multicolumn{2}{c}{$\Pi_{B,0.8}$} & \multicolumn{2}{c}{$\Pi_{B,0.9}$} & \multicolumn{2}{c|}{$\Pi_{B,0.95}$} & \multicolumn{2}{c|}{Ext} & \multicolumn{2}{c}{Dens}\tabularnewline
\hline 
10 & 0&073 & 0&046 & 0&052 & 0&063 & 0&088 & 0&043 & 0&042 & 0&028 & 0&088 & 0&09\tabularnewline
15 & 0&062 & 0&043 & 0&04 & 0&05 & 0&078 & 0&044 & 0&038 & 0&026 & 0&076 & 0&076\tabularnewline
20 & 0&057 & 0&039 & 0&046 & 0&047 & 0&064 & 0&044 & 0&044 & 0&031 & 0&064 & 0&064\tabularnewline
25 & 0&044 & 0&035 & 0&041 & 0&04 & 0&048 & 0&038 & 0&034 & 0&026 & 0&05 & 0&052\tabularnewline
\end{tabular}

\vspace*{.5em}

\begin{tabular}{c|r@{\extracolsep{0pt}.}lr@{\extracolsep{0pt}.}lr@{\extracolsep{0pt}.}lr@{\extracolsep{0pt}.}l|r@{\extracolsep{0pt}.}lr@{\extracolsep{0pt}.}lr@{\extracolsep{0pt}.}lr@{\extracolsep{0pt}.}l|r@{\extracolsep{0pt}.}l|r@{\extracolsep{0pt}.}l}
\multicolumn{21}{l}{\emph{  Unequal mean ($H_1$)}}\tabularnewline
$n_{g}$ & \multicolumn{2}{c}{$L_{A}^{2}$} & \multicolumn{2}{c}{$\Pi_{A,0.8}$} & \multicolumn{2}{c}{$\Pi_{A,0.9}$} & \multicolumn{2}{c|}{$\Pi_{A,0.95}$} & \multicolumn{2}{c}{$L_{B}^{2}$} & \multicolumn{2}{c}{$\Pi_{B,0.8}$} & \multicolumn{2}{c}{$\Pi_{B,0.9}$} & \multicolumn{2}{c|}{$\Pi_{B,0.95}$} & \multicolumn{2}{c|}{Ext} & \multicolumn{2}{c}{Dens}\tabularnewline
\hline 
10 & 0&203 & 0&215 & 0&596 & 0&846 & 0&259 & 0&215 & 0&572 & 0&776 & 0&259 & 0&124\tabularnewline
15 & 0&336 & 0&352 & 0&868 & 0&989 & 0&388 & 0&36 & 0&863 & 0&976 & 0&394 & 0&134\tabularnewline
20 & 0&532 & 0&484 & 0&975 & 0&999 & 0&591 & 0&5 & 0&972 & 0&999 & 0&59 & 0&148\tabularnewline
25 & 0&766 & 0&607 & 0&993 & \multicolumn{2}{c|}{1} & 0&805 & 0&618 & 0&993 & \multicolumn{2}{c|}{1} & 0&794 & 0&173\tabularnewline
\end{tabular}
\end{table}

Under the Equal Mean scenario ($H_0$), the proportion of rejection for all methods were below 0.1 and approached the nominal level $\alpha=0.05$ as $n_g$ increased, indicating that the tests have approximately the correct size. 
The proposed norm-based bootstrap tests $L_{B}^{2}$ was slightly more liberal than its asymptotic version $L_A^2$, while the projection-based $\Pi_{B,0.95}$ was slightly more conservative than $\Pi_{A,0.95}$.
In the Unequal Mean scenario ($H_1$), all tests based on the square root densities, i.e. our proposed intrinsic tests and the extrinsic test, outperformed the bootstrap test based on the original density (last column, \ref{tab:taxi}). 
This highlights the Hilbert sphere $\Sinf$ as a more appropriate geometry for detecting small changes in the demand patterns than a flat space of original densities.
The proposed projection-based tests with FVE threshold $r=0.9$ and $0.95$ had the highest power for all sample sizes, outperforming the intrinsic and extrinsic norm-based bootstrap tests. 

\noindent \textbf{Supplemental Materials}. This manuscript is accompanied by Supplemental Materials including proofs and additional simulations, which will be made available upon request.

\appendix
\section*{Appendix}
\section{Fréchet Derivatives\label{app:fret}} 

The Fr\'{e}chet derivative is reviewed here, following the definitions in Chapter~I of \citet{lang:99}. 
In this section, let $E$, $E_{j},$ $F$ be Banach spaces with norms $\norm{\cdot}_{E}$, $\norm{\cdot}_{E_{j}},$ and $\norm{\cdot}_{F}$, respectively, for $j=1,\dots,p$. 
Denote $\cB(E,F)$ as the space of continuous linear maps from $E$ into $F$, which is a Banach space equipped with the operator norm $\norm g=\sup_{\norm e_{E}=1}\norm{g(e)}_{F}$ for $g\in\cB(E,F)$. 
Also let $\cB(E_{1},\dots,E_{p};F)$ denote the Banach space of multilinear maps equipped with the operator norm
\[
\norm h=\sup_{\norm{e_{1}}_{E_{1}}=\dots=\norm{e_{p}}_{E_{p}}=1}\norm{h(e_{1},\dots,e_{p})}_{F}, 
\]
and write for short $\cB(E^{p},F)=\cB(E,\dots,E;F)$.
Repeated linear operator 
\[
g_{\text{rep}}\in\cB(\underbrace{E,\,\cB(E,\,\dots,\,\cB(E}_{p\text{ times}},F)\dots))
\]
is isometrically identified by a multilinear map $g_{\text{mult}} \in \cB(E^{p},F)$, as 
\[
g_{\text{mult}}(e_{1},\dots,e_{p}) = g_{\text{rep}}(e_{1})\dots(e_{p}).
\]
This identification gives rise to a Banach space isomorphism \citep[Proposition 2.4, p7,][]{lang:99}; we use the same notation to denote both maps. 

Let $f:U\subset E\rightarrow F$ be a continuous map.
\begin{defn}
Function $f:U\subset E\rightarrow F$ is said to be \emph{(Fr\'{e}chet) differentiable}
at a point $x_{0}\in U$ if there exists a continuous linear map $l$
of $E$ into $F$ such that for $y\in E,$
\[
f(x_{0}+y)=f(x_{0})+l(y)+\epsilon(y),
\]
where $\norm{\epsilon(y)}_{F}\tozero$ as $\norm y_{E}\tozero$. 
The linear map $l$ is called the \emph{(Fr\'{e}chet) derivative} of $f$ at $x_{0}$, denoted as $Df(x_0)$. 
If $f$ is differentiable at every point in $U$, then the derivative $Df$ is a map
\[
Df:U\rightarrow\cB(E,F).
\]
\end{defn}
\begin{defn}
Map $f:U\subset E\rightarrow F$ is said to be \emph{directional differentiable} at a point $x_{0}\in U$ if there exists a function $l:E\rightarrow F$ such that 
\[
l(y)=\lim_{t\tozero}\frac{f(x_{0}+ty)-f(x_{0})}{t}
\]
exists for all $y\in E$. 
The linear map $l$ is called the \emph{directional derivative} of $f$ at $x_{0}$.
\end{defn}

If a map $f$ is Fr\'{e}chet differentiable, then it is also directional differentiable and the two derivatives match. 
In what follows, ``differentiability'' refers to Fr\'{e}chet differentiability unless otherwise noted, and the directional differentiation is used for calculation. 
Higher order derivatives and partial derivatives are defined in a recursive fashion.
Since the derivative $Df(x_{0})$ is in $\cB(E,F)$, a Banach space, the \emph{$p$th order derivative} $D^{p}f$ is defined as $D(D^{p-1}f)$, a map of $U$ into $\cB(E,\,\cB(E,\,\dots,\,\cB(E,F)\dots))\simeq\cB(E^{p},F)$.
A map is said to be smooth if the derivatives of all orders exist.
For a bivariate map $h:E_{1}\times E_{2}\rightarrow F$, the \emph{partial
derivative} with respect to the first argument at $(x_{0},y_{0}) \in U \times V \subset E_1 \times E_2$ is denoted as $D_{1}h(x_{0},y_{0}),$ where
\[
D_{1}h:U\times V\rightarrow\cB(E_{1},F),\quad D_{1}h(x_{0},y_{0})=Df_{y_{0}}(x_{0}),
\]
for $f_{y_{0}}(x)=h(x,y_{0})$. 
The partial derivative $D_{2}h$ w.r.t. the second argument is similarly defined.
\begin{prop}[Chain rule, \cite{lang:99}]
\label{prop:chain}If $f:U\rightarrow V$ is differentiable at $x_{0}$,
and $g:V\rightarrow W$ is differentiable at $f(x_{0})$, then $g\circ f$
is differentiable at $x_{0}$, and
\[
D(g\circ f)(x_{0})=Dg\circ f(x_{0})(Df(x_{0})).
\]
\end{prop}
For $f:U\rightarrow\cB(E,F)$ and $g:U\rightarrow\cB(F,G)$ defined
on an open set $U\subset E$, denote $f\cdot g$ as the function $u\mapsto f(u)\circ g(u).$
The chain rule can be compactly written as 
\[
D(g\circ f)=Dg\circ f\cdot Df.
\]

\references


\begin{thebibliography}{50}
\newcommand{\enquote}[1]{``#1''}
\expandafter\ifx\csname natexlab\endcsname\relax\def\natexlab#1{#1}\fi

\bibitem[{Afsari(2011)}]{afsa:11}
Afsari, B. (2011), \enquote{Riemannian
  {{{\emph{L}}}}{\emph{{\textsuperscript{p}}}} Center of Mass: {{Existence}},
  Uniqueness, and Convexity,} \textit{Proceedings of the American Mathematical
  Society}, 139, 655--673.

\bibitem[{{Ahidar-Coutrix} et~al.(2019){Ahidar-Coutrix}, Le~Gouic, and
  Paris}]{ahid:19}
{Ahidar-Coutrix}, A., Le~Gouic, T., and Paris, Q. (2019), \enquote{Convergence
  Rates for Empirical Barycenters in Metric Spaces: Curvature, Convexity and
  Extendable Geodesics,} \textit{Probability Theory and Related Fields}.

\bibitem[{Aue et~al.(2018)Aue, Rice, and S{\"o}nmez}]{aue:18}
Aue, A., Rice, G., and S{\"o}nmez, O. (2018), \enquote{Detecting and Dating
  Structural Breaks in Functional Data without Dimension Reduction,}
  \textit{Journal of the Royal Statistical Society: Series B (Statistical
  Methodology)}, 80, 509--529.

\bibitem[{Bauer et~al.(2017)Bauer, Eslitzbichler, and Grasmair}]{baue:17}
Bauer, M., Eslitzbichler, M., and Grasmair, M. (2017), \enquote{Landmark-Guided
  Elastic Shape Analysis of Human Character Motions,} \textit{Inverse Problems
  \& Imaging}, 11, 601--621.

\bibitem[{Berkes et~al.(2009)Berkes, Gabrys, Horv{\'a}th, and
  Kokoszka}]{berk:09}
Berkes, I., Gabrys, R., Horv{\'a}th, L., and Kokoszka, P. (2009),
  \enquote{Detecting Changes in the Mean of Functional Observations,}
  \textit{Journal of the Royal Statistical Society: Series B (Statistical
  Methodology)}, 71, 927--946.

\bibitem[{Bhattacharya and Lin(2017)}]{bhat:17}
Bhattacharya, R. and Lin, L. (2017), \enquote{Omnibus {{CLTs}} for
  {{Fr\'echet}} Means and Nonparametric Inference on Non-{{Euclidean}} Spaces,}
  \textit{Proceedings of the American Mathematical Society}, 145, 413--428.

\bibitem[{Bhattacharya and Patrangenaru(2003)}]{bhat:03}
Bhattacharya, R. and Patrangenaru, V. (2003), \enquote{Large Sample Theory of
  Intrinsic and Extrinsic Sample Means on Manifolds - {{I}},} \textit{Annals of
  Statistics}, 31, 1--29.

\bibitem[{Bhattacharya and Patrangenaru(2005)}]{bhat:05}
--- (2005), \enquote{Large Sample Theory of Intrinsic and Extrinsic Sample
  Means on Manifolds - {{II}},} \textit{Annals of statistics}, 33, 1225--1259.

\bibitem[{Cheng et~al.(2009)Cheng, Ghosh, Jiang, and Deriche}]{chen:09-1}
Cheng, J., Ghosh, A., Jiang, T., and Deriche, R. (2009), \enquote{A
  {{Riemannian Framework}} for {{Orientation Distribution Function
  Computing}},} in \textit{Medical {{Image Computing}} and
  {{Computer}}-{{Assisted Intervention}} \textendash{} {{MICCAI}} 2009}, eds.
  Yang, G.-Z., Hawkes, D., Rueckert, D., Noble, A., and Taylor, C., {Berlin,
  Heidelberg}: {Springer Berlin Heidelberg}, vol. 5761, pp. 911--918.

\bibitem[{Chu and Chen(2019)}]{chu:19}
Chu, L. and Chen, H. (2019), \enquote{Asymptotic Distribution-Free Change-Point
  Detection for Multivariate and Non-{{Euclidean}} Data,} \textit{The Annals of
  Statistics}, 47, 382--414.

\bibitem[{Dai and M{\"u}ller(2018)}]{dai:18-5}
Dai, X. and M{\"u}ller, H.-G. (2018), \enquote{Principal Component Analysis for
  Functional Data on {{Riemannian}} Manifolds and Spheres,} \textit{Annals of
  Statistics}, 46, 3334--3361.

\bibitem[{Du et~al.(2014)Du, Goh, Kushnarev, and Qiu}]{du:14}
Du, J., Goh, A., Kushnarev, S., and Qiu, A. (2014), \enquote{Geodesic
  Regression on Orientation Distribution Functions with Its Application to an
  Aging Study,} \textit{NeuroImage}, 87, 416--426.

\bibitem[{Dubey and M{\"u}ller(2019)}]{dube:19-1}
Dubey, P. and M{\"u}ller, H.-G. (2019), \enquote{Fr\'echet Analysis of Variance
  for Random Objects,} \textit{Biometrika}, 106, 803--821.

\bibitem[{Dubey and M{\"u}ller(2020)}]{dube:20}
--- (2020), \enquote{Functional Models for Time-Varying Random Objects,}
  \textit{Journal of the Royal Statistical Society: Series B (Statistical
  Methodology)}, 82, 275--327.

\bibitem[{Ellingson et~al.(2013)Ellingson, Patrangenaru, and
  Ruymgaart}]{elli:13}
Ellingson, L., Patrangenaru, V., and Ruymgaart, F. (2013),
  \enquote{Nonparametric Estimation of Means on {{Hilbert}} Manifolds and
  Extrinsic Analysis of Mean Shapes of Contours,} \textit{Journal of
  Multivariate Analysis}, 122, 317--333.

\bibitem[{Eltzner and Huckemann(2019)}]{eltz:19}
Eltzner, B. and Huckemann, S.~F. (2019), \enquote{A Smeary Central Limit
  Theorem for Manifolds with Application to High-Dimensional Spheres,}
  \textit{The Annals of Statistics}, 47, 3360--3381.

\bibitem[{Fletcher et~al.(2004)Fletcher, Lu, Pizer, and Joshi}]{flet:04}
Fletcher, P.~T., Lu, C., Pizer, S.~M., and Joshi, S. (2004), \enquote{Principal
  Geodesic Analysis for the Study of Nonlinear Statistics of Shape,}
  \textit{IEEE Transactions on Medical Imaging}, 23, 995--1005.

\bibitem[{Fr{\'e}chet(1948)}]{frec:48}
Fr{\'e}chet, M. (1948), \enquote{Les \'El\'ements Al\'eatoires de Nature
  Quelconque Dans Un Espace Distanci\'e,} \textit{Annales de l'Institut Henri
  Poincar\'e}, 10, 215--310.

\bibitem[{Gouic et~al.(2019)Gouic, Paris, Rigollet, and Stromme}]{goui:19}
Gouic, T.~L., Paris, Q., Rigollet, P., and Stromme, A.~J. (2019), \enquote{Fast
  Convergence of Empirical Barycenters in {{Alexandrov}} Spaces and the
  {{Wasserstein}} Space,} \textit{arXiv:1908.00828 [math, stat]}.

\bibitem[{H{\'a}jek(1962)}]{haje:62}
H{\'a}jek, J. (1962), \enquote{On Linear Statistical Problems in Stochastic
  Processes,} \textit{Czechoslovak Mathematical Journal}, 12, 404--444.

\bibitem[{Hall(1992)}]{hall:92-2}
Hall, P. (1992), \textit{The {{Bootstrap}} and {{Edgeworth Expansion}}}, {New
  York}: {Springer}.

\bibitem[{Henning and Srivastava(2016)}]{henn:16}
Henning, W. and Srivastava, A. (2016), \enquote{A Two-Sample Test for
  Statistical Comparisons of Shape Populations,} in \textit{2016 {{IEEE Winter
  Conference}} on {{Applications}} of {{Computer Vision}} ({{WACV}})}, {IEEE},
  pp. 1--9.

\bibitem[{Horv{\'a}th et~al.(2013)Horv{\'a}th, Kokoszka, and
  Reeder}]{horv:13-1}
Horv{\'a}th, L., Kokoszka, P., and Reeder, R. (2013), \enquote{Estimation of
  the Mean of Functional Time Series and a Two-Sample Problem,} \textit{Journal
  of the Royal Statistical Society: Series B (Statistical Methodology)}, 75,
  103--122.

\bibitem[{Hotz and Huckemann(2015)}]{hotz:15}
Hotz, T. and Huckemann, S. (2015), \enquote{Intrinsic Means on the Circle:
  {{Uniqueness}}, Locus and Asymptotics,} \textit{Annals of the Institute of
  Statistical Mathematics}, 67, 177--193.

\bibitem[{Hsing and Eubank(2015)}]{hsin:15}
Hsing, T. and Eubank, R. (2015), \textit{Theoretical {{Foundations}} of
  {{Functional Data Analysis}}, with an {{Introduction}} to {{Linear
  Operators}}}, {Hoboken}: {Wiley}.

\bibitem[{Huckemann(2012)}]{huck:12}
Huckemann, S.~F. (2012), \enquote{On the Meaning of Mean Shape: Manifold
  Stability, Locus and the Two Sample Test,} \textit{Annals of the Institute of
  Statistical Mathematics}, 64, 1227--1259.

\bibitem[{Joshi et~al.(2007)Joshi, Klassen, Srivastava, and Jermyn}]{josh:07}
Joshi, S.~H., Klassen, E., Srivastava, A., and Jermyn, I. (2007), \enquote{A
  Novel Representation for {{Riemannian}} Analysis of Elastic Curves in
  {{{\emph{R}}}}{\emph{{\textsuperscript{n}}}},} in \textit{2007 {{IEEE
  Conference}} on {{Computer Vision}} and {{Pattern Recognition}}}, {IEEE}, pp.
  1--7.

\bibitem[{Karcher(1977)}]{karc:77}
Karcher, H. (1977), \enquote{Riemannian Center of Mass and Mollifier
  Smoothing,} \textit{Communications on Pure and Applied Mathematics}, 30,
  509--541.

\bibitem[{Lang(1999)}]{lang:99}
Lang, S. (1999), \textit{Fundamentals of Differential Geometry}, {New York}:
  {Springer}.

\bibitem[{Lazar and Lin(2017)}]{laza:17}
Lazar, D. and Lin, L. (2017), \enquote{Scale and Curvature Effects in Principal
  Geodesic Analysis,} \textit{Journal of Multivariate Analysis}, 153, 64--82.

\bibitem[{Le(2001)}]{le:01}
Le, H. (2001), \enquote{Locating {{Fr\'echet}} Means with Application to Shape
  Spaces,} \textit{Advances in Applied Probability}, 33, 324--338.

\bibitem[{Lin and Yao(2019)}]{lin:19}
Lin, Z. and Yao, F. (2019), \enquote{Intrinsic {{Riemannian}} Functional Data
  Analysis,} \textit{The Annals of Statistics}, 47, 3533--3577.

\bibitem[{Petersen et~al.(2019)Petersen, Liu, and Divani}]{pete:19-1}
Petersen, A., Liu, X., and Divani, A.~A. (2019), \enquote{Wasserstein
  {{{\emph{F}}}}-Tests and {{Confidence Bands}} for the {{Fr\'echet
  Regression}} of {{Density Response Curves}},} \textit{arXiv}.

\bibitem[{Petersen and M{\"u}ller(2016)}]{pete:16}
Petersen, A. and M{\"u}ller, H.-G. (2016), \enquote{Functional Data Analysis
  for Density Functions by Transformation to a {{Hilbert}} Space,} \textit{The
  Annals of Statistics}, 44, 183--218.

\bibitem[{Rao(1945)}]{rao:45}
Rao, C.~R. (1945), \enquote{Information and the Accuracy Attainable in the
  Estimation of Statistical Parameters,} \textit{Bulletin of Calcutta
  Mathematical Society}, 81--91.

\bibitem[{Rudin(1973)}]{rudi:73}
Rudin, W. (1973), \textit{Functional {{Analysis}}}, {New York}: {McGraw-Hill}.

\bibitem[{Sch{\"o}tz(2019)}]{scho:19}
Sch{\"o}tz, C. (2019), \enquote{Convergence Rates for the Generalized
  {{Fr\'echet}} Mean via the Quadruple Inequality,} \textit{Electronic Journal
  of Statistics}, 13, 4280--4345.

\bibitem[{Srivastava et~al.(2007)Srivastava, Jermyn, and Joshi}]{sriv:07}
Srivastava, A., Jermyn, I., and Joshi, S. (2007), \enquote{Riemannian Analysis
  of Probability Density Functions with Applications in Vision,} in
  \textit{2007 {{IEEE Conference}} on {{Computer Vision}} and {{Pattern
  Recognition}}}, pp. 1--8.

\bibitem[{Strait et~al.(2019)Strait, Chkrebtii, and Kurtek}]{stra:19}
Strait, J., Chkrebtii, O., and Kurtek, S. (2019), \enquote{Automatic Detection
  and Uncertainty Quantification of Landmarks on Elastic Curves,}
  \textit{Journal of the American Statistical Association}, 1--23.

\bibitem[{Su et~al.(2014)Su, Kurtek, Klassen, and Srivastava}]{su:14}
Su, J., Kurtek, S., Klassen, E., and Srivastava, A. (2014),
  \enquote{Statistical Analysis of Trajectories on {{Riemannian}} Manifolds:
  {{Bird}} Migration, Hurricane Tracking and Video Surveillance,} \textit{The
  Annals of Applied Statistics}, 8, 530--552.

\bibitem[{Tucker et~al.(2013)Tucker, Wu, and Srivastava}]{tuck:13}
Tucker, J.~D., Wu, W., and Srivastava, A. (2013), \enquote{Generative Models
  for Functional Data Using Phase and Amplitude Separation,}
  \textit{Computational Statistics \& Data Analysis}, 61, 50--66.

\bibitem[{{van der Vaart} and Wellner(1996)}]{van:96}
{van der Vaart}, A. and Wellner, J. (1996), \textit{Weak {{Convergence}} and
  {{Empirical Processes}}: {{With Applications}} to {{Statistics}}}, {New
  York}: {Springer}.

\bibitem[{Wang et~al.(2016)Wang, Chiou, and M{\"u}ller}]{wang:16}
Wang, J.-L., Chiou, J.-M., and M{\"u}ller, H.-G. (2016), \enquote{Functional
  Data Analysis,} \textit{Annual Review of Statistics and its Application}, 3,
  257--295.

\bibitem[{Wellner and Zhan(1996)}]{well:96}
Wellner, J.~A. and Zhan, Y. (1996), \enquote{Bootstrapping
  {{{\emph{Z}}}}-Estimators,} Tech. rep.

\bibitem[{Wu and Srivastava(2014)}]{wu:14}
Wu, W. and Srivastava, A. (2014), \enquote{Analysis of Spike Train Data:
  {{Alignment}} and Comparisons Using the Extended {{Fisher}}-{{Rao}} Metric,}
  \textit{Electronic Journal of Statistics}, 8, 1776--1785.

\bibitem[{Xie et~al.(2017)Xie, Kurtek, Bharath, and Sun}]{xie:17}
Xie, W., Kurtek, S., Bharath, K., and Sun, Y. (2017), \enquote{A Geometric
  Approach to Visualization of Variability in Functional Data,} \textit{Journal
  of the American Statistical Association}, 112, 979--993.

\bibitem[{Younes(1998)}]{youn:98}
Younes, L. (1998), \enquote{Computable Elastic Distances between Shapes,}
  \textit{SIAM Journal on Applied Mathematics}, 58, 565--586.

\bibitem[{Yu et~al.(2017)Yu, Lu, and Marron}]{yu:17}
Yu, Q., Lu, X., and Marron, J. (2017), \enquote{Principal Nested Spheres for
  Time-Warped Functional Data Analysis,} \textit{Journal of Computational and
  Graphical Statistics}, 26, 144--151.

\bibitem[{Zhu et~al.(2009)Zhu, Chen, Ibrahim, Li, Hall, and Lin}]{zhu:09}
Zhu, H., Chen, Y., Ibrahim, J.~G., Li, Y., Hall, C., and Lin, W. (2009),
  \enquote{Intrinsic Regression Models for Positive-Definite Matrices with
  Applications to Diffusion Tensor Imaging,} \textit{Journal of the American
  Statistical Association}, 104, 1203--1212.

\bibitem[{Ziezold(1977)}]{ziez:77}
Ziezold, H. (1977), \enquote{On Expected Figures and a Strong Law of Large
  Numbers for Random Elements in Quasi-Metric Spaces,} in \textit{Transactions
  of the {{Seventh Prague Conference}} on {{Information Theory}}, {{Statistical
  Decision Functions}}, {{Random Processes}} and of the 1974 {{European
  Meeting}} of {{Statisticians}}}, pp. 591--602.

\end{thebibliography}
\end{document}